\documentclass[12pt,twoside,final,psamsfonts]{amsart}

\usepackage[psamsfonts]{amssymb}
\usepackage{times,a4wide}
\theoremstyle{plain}         
\newtheorem*{theorem*}{Theorem}

\newtheorem{theorem}{Theorem}[section]
\newtheorem{proposition}[theorem]{Proposition}
\newtheorem*{proposition*}{Proposition}
\newtheorem{corollary}[theorem]{Corollary}
\newtheorem*{corollary*}{Corollary}

\newtheorem*{lemma*}{Lemma}

\theoremstyle{definition}
\newtheorem{remark}[theorem]{Remark}
\newtheorem*{remark*}{Remark}
\newtheorem{example}[theorem]{Example}

\theoremstyle{definition}

\newtheorem*{definition*}{Definition}
\newcommand{\nc}{\newcommand}

\newcommand{\N}{{\mathbb N}}
\newcommand{\D}{{\mathbb D}}
\newcommand{\R}{{\mathbb R}}

\newcommand{\T}{{\mathbb T}}

\newcommand{\Int}{\operatorname{Int}}
\newcommand{\Hol}{\operatorname{Hol}}

\newcommand{\lra}{\longrightarrow}

\newcommand{\eps}{\varepsilon}
\newcommand{\vp}{\varphi}
\nc{\bea}{\begin{eqnarray}} 
\nc{\eea}{\end{eqnarray}}
\nc{\beqa}{\begin{eqnarray*}} 
\nc{\eeqa}{\end{eqnarray*}}
\nc{\Hi}{H^{\infty}}
\nc{\loi}{\ell^{\infty}}
\nc{\NL}{N^+\vert \Lambda}
\nc{\hf}{{\mathcal H}_{\vp}}
\nc{\liL}{\lambda\in\Lambda}
\nc{\nn}{\nonumber}

\renewcommand{\qedsymbol}{$\blacksquare$}

\title[Interpolating sequences for the Nevanlinna and Smirnov classes]
{Interpolating sequences for the Nevanlinna and Smirnov classes}

\author[A. Hartmann, X. Massaneda, A. Nicolau]{Andreas Hartmann, 
Xavier Massaneda \& Artur Nicolau}

\address{Laboratoire de Math\'ematiques Pures de Bordeaux,
Universit\'e Bordeaux I, 351 cours de la Lib\'eration,
33405 Talence, France}

\address{Departament de Matem\`atica Aplicada i An\`alisi,
Universitat  de Barcelona, Gran Via 585, 08071-Bar\-ce\-lo\-na, Spain}

\address{Departament de Matem\`atiques, Facultat de Ci\`encies,
Universitat Aut\`onoma de Barcelona, 08193-Bellaterra, Spain}

\thanks{All authors supported by a PICS program of Generalitat de Catalunya 
and CNRS. First and third author also supported  by European Commission
Research Training Network HPRN-CT-2000-00116.
Second author also supported by the DGICYT grant BFM2002-04072-C02-01 and the 
CIRIT grant 2001-SGR00172. Third author also supported by DGICYT grant 
BFM2002-00571 and the CIRIT grant 2001-SGR00431}
\email{hartmann@math.u-bordeaux.fr, xavier@mat.ub.es, artur@mat.uab.es}
\date{Decembre 5, 2002}



\begin{document}

\begin{abstract} We give analytic Carleson-type
characterisations of the interpolating
sequences for the Nevanlinna and Smirnov classes. From this we deduce 
necessary and  sufficient geometric conditions, both expressed in terms of a
certain non-tangential maximal function associated to the sequence. Some
examples show that the gap between the necessary and the sufficient conditions
cannot be covered. We also discuss the relationship between our results and
the previous work of Naftalevi\v c for the Nevanlinna class, and Yanagihara
for the Smirnov class. Finally, we observe that the arguments used in the
previous proofs show that interpolating sequences for ``big'' Hardy-Orlicz
spaces are in general different from  those for the scale included in the
classical Hardy spaces.   
\end{abstract}

\maketitle

\section{Introduction and basic definitions}\label{intro}

Let $\Lambda$ be a discrete sequence of points in the unit disk $\D$. For a 
space of holomorphic functions $X$, the interpolation problem 
consists of describing the trace of $X$ on $\Lambda$, i.e.\ the set of
restrictions  $X\vert \Lambda$, which is regarded as a sequence space. One way
of defining  interpolating sequences is to fix a priori a natural trace space $l$
and look  for conditions ensuring that $X\vert \Lambda=l$. The second
possibility is  to require the trace space 
\[
X\vert \Lambda=\{(f(\lambda))_{\lambda\in\Lambda} : f\in X\}
\] 
to be {\it ideal}, i.e. $\ell^\infty  X\vert \Lambda\subset X\vert \Lambda$
(see definitions below). This approach is motivated by the property of
unconditional bases to be absolutely convergent (see \cite[Section C.3.1 (Volume 2)]{Nik02} 
for more about this, in particular, Theorem C.3.1.4) 
and it is natural at least for those spaces that are stable under
multiplication by $\Hi$, the space of bounded holomorphic functions
on $\D$. 

For many spaces (for instance Hardy and Bergman spaces), both definitions turn
out to be equivalent, provided that the a priori fixed trace space is chosen in
a natural way. The situation changes for the non-Banach classes we have in
mind. In order to illustrate this, we briefly discuss the known results
for the spaces we will deal with, the \emph{Nevanlinna class}
\[
N=\bigl\{f\in \Hol(\D):\lim_{r \to 1}\frac{1}{2\pi}\int_{-\pi}^{\pi}
 \log^+|f(re^{it})|\;dt<\infty\bigr\}
\]
and the related \emph{Smirnov class }
\[
 N^+=\bigl\{f\in N:\lim_{r \to 1}\frac{1}{2\pi}\int_{-\pi}^{\pi}
 \log^+|f(re^{it})|\;dt=\frac{1}{2\pi}\int_{-\pi}^{\pi}
 \log^+|f(e^{it})|\;dt\bigr\}.
\]

In 1956, in a very interesting article, Naftalevi\v c \cite{Na56} described the
sequences $\Lambda$ for which the trace $N\vert \Lambda$ coincides with the
sequence space $l_{\text{Na}}=\{(a_\lambda)_\lambda: \sup_{\lambda}
(1-|\lambda|)\log^+|a_\lambda|<\infty\}$. The choice of $l_{\text{Na}}$ is 
motivated by the fact that $\sup_z (1-|z|)\log^+ |f(z)|<\infty$ for $f\in N$,
and this growth is attained. 
Unfortunately, the growth condition imposed in $l_{\text{Na}}$  forces the 
sequences to be confined in a finite union of Stolz angles.  Consequently a big
class of Carleson sequences (i.e.\ sequences such that 
$\Hi|\Lambda=l^{\infty}$), 
namely
those tending tangentially to the boundary,  cannot be interpolating in the sense
of Naftalevi\v c. This does not seem  natural, for $\Hi$ is in the multiplier
space of $N$ (since $N$ is an  algebra, its multiplier space  obviously coincides
with itself).  Further comments on Naftalevi\v c's result can be found in
\cite{HM2}. 

A similar problem occurs in the Smirnov class. In \cite{yana2}, Yanagihara proved
that in order that $N^+\vert \Lambda$  contains the space
$l_{\text{Ya}}=\{(a_\lambda)_\lambda:\sum_\lambda (1-|\lambda|)
\log^+|a_\lambda|<\infty\}$,  it is sufficient that $\Lambda$ is a Carleson
sequence (he also gave a  necessary condition that we will discuss below).
However there are Carleson sequences such that $N^+\vert \Lambda$ does not embed
into  $l_{\text{Ya}}$ \cite[Theorem 3]{yana2}. 

In conclusion, it seems quite difficult to find a ``natural'' trace for these
spaces. Therefore we consider the following definitions.

\begin{definition*}\label{ideal}
A sequence space $l$ is called {\it ideal} if $\ell^\infty l\subset l$, i.e.
whenever $(a_n)_n\in l$ and  $ (\omega_n)_n\in \ell^\infty$,  then also $(\omega_n
a_n)_n\in l$.
\end{definition*}

\begin{definition*}\label{freeint}
Let $X$ be a space of holomorphic functions in $\D$. 
A sequence $\Lambda\subset \D$ is called {\it
free interpolating for $X$} if $X\vert \Lambda$ is ideal.
We denote $\Lambda\in \Int X$.
\end{definition*}

Since the Nevanlinna and Smirnov classes contain the constants,  free
interpolation for these classes entails the existence of a nonzero function 
$f\in N$ vanishing on $\Lambda$, hence the Blaschke condition $\sum_\lambda
(1-|\lambda|)<\infty$ is necessary and will be assumed throughout this
paper.

\begin{remark}\label{l-infinit}
For any function algebra $X$ containing the constants,
$X\vert\Lambda$ is ideal if and only if
\beqa 
 \loi\subset X\vert\Lambda.
\eeqa

The inclusion is obviously necessary. In order to see that it is sufficient
notice that, by assumption, for any $(\omega_\lambda)_\lambda\in \loi$  there
exists $g\in X$ such that $g(\lambda)=\omega_\lambda$. Thus, if
$(f(\lambda))_\lambda\in X\vert \Lambda$, the sequence of values
$(\omega_\lambda f(\lambda))_\lambda$ can be interpolated by  $fg\in X$. 
\end{remark}

It is then clear that $\Int N^+\subset\Int N$.

We shall use Remark \ref{l-infinit} not only for the classes $N$ and $N^+$, but
also for ``big'' Hardy-Orlicz classes (see Section \ref{nevanlinna}).

In order to state our results we need to recall some standard facts about the
structure of the Nevanlinna and Smirnov classes (we refer to \cite{Gar},
\cite{Nik02} or \cite{RR} for general refe\-rences). Let $d\sigma$ denote the
normalised Lebesgue measure in the unit circle $\mathbb T$.

A function $f$ is called \emph{outer} if it can be written in the form
\[
 f(z)=C \exp \left\{\int_\T \frac{\zeta+z}{\zeta-z} 
 \log v(\zeta) d\sigma(\zeta) \right\},
\] 
where $|C|=1$, $v>0$  a.e.\ on $\T$ and $\log v \in L^1(\T)$. Such a function is
the quotient $f=f_1/f_2$ of two bounded outer functions  $f_1,f_2\in \Hi$
with  $\|f_i\|_{\infty}\le 1$, $i=1,2$. In particular, the weight $v$ is given by
the boundary values of $|f_1/f_2|$. Setting $w=\log v$, we have
\bea\label{outerfct} 
\log |f(z)|=P[w](z):=\int_\T P(z,\zeta) w(\zeta)d\sigma(\zeta), 
\eea
where 
\[
P(z,\zeta)=\textrm{Re}\left(\frac{\zeta+z}{\zeta-z}\right)=\frac{1-|z|^2}{|\zeta-z|^2}
\] 
is the Poisson kernel in $\D$. Formula (\ref{outerfct}) 
allows to freely switch between 
assertions on outer  function $f$ and assertions on the associated measures 
$w d\sigma$.

Another important family of functions in this context are \emph{inner} functions: 
$I\in \Hi$ such that $|I|=1$ almost everywhere on $\T$. 
Any inner function $I$ can be factorised into
a Blaschke product $B_\Lambda= \prod_n b_{\lambda_n}$ carrying 
the zeros $\Lambda=\{\lambda_n\}_n$  of $I$, where $b_{\lambda}(z)=
\frac{|\lambda|}{\lambda}\frac{\lambda-z}{1-\overline{\lambda}z}$ denotes the usual
M\"obius transformation, and a singular inner function $S$ defined by
\beqa
 S(z)=\exp\left\{-\int_\T \frac{\zeta+z}{\zeta-z}\,d\mu(\zeta)\right\}, 
\eeqa
for some positive Borel measure $\mu$ singular with respect to Lebesgue
measure. 

According to the Riesz-Smirnov factorisation, any
function $f\in N^+$ is represented as
\[
 f=\alpha \frac{B S f_1}{f_2},
\]
where $f_1, f_2$ are outer with $\|f_1\|_{\infty},\|f_2\|_{\infty}\le 1$, 
$S$ is singular inner, $B$ is a Blaschke product and $|\alpha|=1$. Similarly,
functions $f\in N$ are represented as
\[
 f=\alpha \frac{B S_1 f_1}{S_2 f_2},
\]
with $f_i$ outer, $\|f_i\|_{\infty}\le 1$, $S_i$ singular inner,
$B$ is a Blaschke product and $|\alpha|=1$.

Given the Blaschke product $B$ with zero-sequence $\Lambda$, denote
$B_\lambda=B/b_\lambda$. 

Our main results are the following. 

\begin{theorem}\label{thmCNS} 
Let $\Lambda$ be a sequence in $\D$. The
following statements are equivalent:
\begin{itemize}
\item[(a)] $\Lambda\in\Int N^+$.

\item[(b)] There exists an outer function $g$ such that 
\bea\label{CNS}
 |B_\lambda(\lambda)|\geq |g(\lambda)|\quad
 \text{ for all }\lambda\in\Lambda.
\eea

\item[(c)] There exists a positive function $w\in L^1(\T)$ such that
\bea\label{CNSmu}
 \log\frac 1{|B_\lambda(\lambda)|}\leq P[w](\lambda) \quad\text{ for all }
 \lambda\in\Lambda .
\eea
\item[(d)] The trace space is given by:
\beqa
 N^+|\Lambda=l_{N^+}:=\{(a_{\lambda})_{\lambda}&:&\text{ there exists a function 
$u\in L^1(\T)$ such that }\\
 & &P[u](\lambda)\ge \log^+|a_{\lambda}|,\ \liL \}.
\eeqa
\end{itemize}
\end{theorem}

Note that in (b) there is no harm in assuming that $g\in\Hi$ with $\|g\|_{\infty}\le 1$,
i.e. $w<0$. This is easily achieved by neglecting the denominator in the quotient
$g=g_1/g_2$, $\|g_i\|_\infty\leq 1$. Then it is clear from \eqref{outerfct} above that
the assertions (b) and (c) are essentially the same. 

In case $|g|$ is uniformly bounded from below by a positive constant
on the sequence $\Lambda$, condition \eqref{CNS}
is nothing but the classical Carleson condition for interpolation in 
$\Hi$.  

According to the Riesz-Smirnov factorisation described above,  the essential
difference bet\-ween Nevanlinna and Smirnov functions is the extra singular factor
appearing in the denominator in the Nevanlinna case. This is reflected in the
corresponding result for free interpolation in $N$.

Given a measure $\mu$ on $\T$, let $P[\mu](z)$ denote the integral obtained from
\eqref{outerfct} by replacing $w\, d\sigma$ by $d\mu$.

\begin{theorem}\label{thmCNSN}
Let $\Lambda$ be a sequence in $\D$. The following statements are 
equivalent:
\begin{itemize}
\item[(a)] $\Lambda\in\Int N$.

\item[(b)] There exist an outer function $g$ and a singular inner function
$S$ such that 
\bea\label{CNSN}
 |B_\lambda(\lambda)|\geq |g(\lambda)S(\lambda)|\quad
 \text{ for all }\lambda\in\Lambda.
\eea

\item[(c)] There exists a positive finite measure $\mu$ on $\T$ such that
\bea\label{CNSNmu}
 \log\frac 1{|B_\lambda(\lambda)|}\leq P[\mu](\lambda)\quad\text{ for all }
 \lambda\in\Lambda .
\eea
\item[(d)] The trace space is given by:
\beqa
 N|\Lambda=l_N:=\{(a_{\lambda})_{\lambda}&:&\text{ there exists a positive finite
 measure $\mu$ on $\T$ such that }\\
 & &P[\mu](\lambda)\ge \log^+|a_{\lambda}|,\ \liL \}.
\eeqa
\end{itemize}
\end{theorem}

As before, in (b) we can assume $g=g_1\in H^\infty$, $\|g_1\|_{\infty}\le 1$. Let $\mu_S$
denote the singular measure associated with $S$, then the equivalence of (b)  and (c)
can be checked using the measure $d\mu=\log (1/|g_1|)\,dm+d\mu_S$.

By Smirnov's theorem, any positive harmonic function on $\D$ is the Poisson extension
of some positive finite measure on $\T$. Thus our results are closely related
to the following problem.

\textit{Problem}. Given a Blaschke sequence $\Lambda$, describe the
sequences of positive values $\{m_\lambda\}_\lambda$ for which there exists a
positive harmonic function $u$ in the unit disk such that $u(\lambda)\geq
m_\lambda$ for all $\lambda\in\Lambda$. 

A geometric description of such sequences must depend on both the behaviour of
$\{m_\lambda\}_\lambda$ and the geometry of $\Lambda$. We only have an answer
in two extremal cases: when $\Lambda$ approaches the unit circle
non-tangentially,  and when it does it very tangentially (in the sense that
the arcs $\{\zeta\in\T : |\arg \zeta -\arg\lambda|<(1-|\lambda|)^{1/2}\}$ are
pairwise disjoint).

In what follows we would like to see some geometric implications of the 
analytic conditions above. To begin with, we would like to state the maybe surprising
result that separated sequences are interpolating for the Smirnov class (and hence
the Nevanlinna class).
Recall that a sequence $\Lambda$ is called \emph{separated} if
\[
\delta(\Lambda):=\inf_{\lambda\neq \lambda'}|b_{\lambda}(\lambda')| >0. 
\]
For such sequences there always
exists an outer function satisfying \eqref{CNS}
(see Proposition \ref{propsep}),  
thus the following corollary is
immediate from Theorem~\ref{thmCNS}.

\begin{corollary}\label{separated}
Let $\Lambda$ be a separated Blaschke sequence. Then $\Lambda\in\Int N^+$
(and hence $\Lambda\in N$).
\end{corollary}

More precise geometric conditions can be given in terms 
of a non-tangential 
maximal function associated with $\Lambda$.
For $\zeta\in \T$ and $\alpha>1$ define the Stolz angle 
\[
\Gamma_{\alpha}(\zeta)=
\{z\in\D: |z-\zeta|\le \alpha (1-|z|^2)\}.
\] 
In our considerations the
value of $\alpha$ is of no importance, so we will write $\Gamma(\zeta)$ 
for the generic Stolz angle with fixed aperture $\alpha$. 
For a given $\Lambda$
consider now the non-tangential maximal function  
\[
 M_\Lambda(\zeta)=\sup_{\lambda\in \Gamma(\zeta)} 
 \log\frac 1{|B_\lambda(\lambda)|}.
\]

Let $L_w^{1}$ denote the weak-$L^1$ space and let
\[
 L^{1}_{w,0}(\T)=\{ f: \lim_{t\to \infty}t \sigma(\{\zeta : |f(\zeta)|>t\}|)=0\}. 
\]

\begin{corollary}\label{maxfct}
Let $\Lambda $ be a sequence in $\D$. 

\begin{itemize}
\item[(a)] If $\Lambda\in\Int N^+$ then $M_\Lambda\in L^1_{w,0}(\T)$.
If $\Lambda\in\Int N$ then $M_\Lambda\in L^1_{w}(\T)$.
\item[(b)] If $M_\Lambda\in L^1(\T)$ then $\Lambda\in\Int N^+$ (and hence
$\Lambda\in \Int N$).
\end{itemize}

\end{corollary}

\begin{proof} 
(a) is a consequence of (c) in Theorem \ref{thmCNSN}. Indeed, it is a general 
fact that the non-tangential maximal function of the Poisson transform of  a
positive finite measure belongs to $L^{1}_{w}(\T)$ 
(see for instance \cite[p.28-29]{Gar}).  A
more careful analysis of the cited result shows that if $\mu$ is  absolutely
continuous, then its Poisson transform is in $L^{1}_{w,0}(\T)$, which yields
the result for $N^+$.

In order to prove (b) consider the arcs associated with $\lambda\in\D$, 
defined as
\begin{equation}\label{shadow}
I_{\lambda}=\{\zeta\in\T:|\arg\zeta-\arg\lambda|\le \pi(1-|\lambda|)\}
\end{equation}
Take $w=(1+\pi^2) M_\Lambda$ and apply (c) in Theorem~\ref{thmCNS}:
\begin{eqnarray}\label{poissonestimate}
  P[ w](\lambda)\geq
 \frac 1{1-|\lambda|} \int_{I_\lambda} M_\Lambda(\zeta) d\sigma(\zeta)\nn 
 \geq \frac 1{1-|\lambda|} \int_{I_\lambda}  \log\frac
 1{|B_\lambda(\lambda)|} d\sigma(\zeta)=\log\frac
 1{|B_\lambda(\lambda)|}.
\end{eqnarray}
\end{proof}

Some Carleson-type conditions can be deduced from the implicit analytic
characterisation and from Corollary \ref{maxfct}. 

\begin{corollary}\label{CNgeom}
({\rm a}) If $\Lambda\in\Int N^+$, then 
\bea\label{CN}
 \lim_{|\lambda|\to 1} (1-|\lambda|)\log\frac 1{  |B_{\lambda}(\lambda)|}=0.
\eea

({\rm b}) If $\Lambda\in\Int N$, then
\bea\label{CNN}
 \sup_{\lambda\in\Lambda} (1-|\lambda|)\log\frac 1{  |B_{\lambda}(\lambda)|}<\infty.
\eea
\end{corollary}

\begin{proof}
Since
\[
 I_\lambda\subset \bigl\{\zeta\in\T : M_\Lambda(\zeta)\geq \log\frac
 1{|B_\lambda(\lambda)|}\bigr\} ,\quad \lambda\in\Lambda,
\]
it suffices to apply
condition (a) of Corollary \ref{maxfct}.
\end{proof}

These are the best possible necessary conditions
expressed in these terms, as reveals the next result.

\begin{proposition}\label{thmnafta}
Assume that $\Lambda\subset\D$ lies in a
finite union of Stolz angles. 
\begin{itemize}
\item[(a)] $\Lambda\in\Int N^+$ if and only if \eqref{CN} holds.

\item[(b)] $\Lambda\in\Int N$ if and only if \eqref{CNN} holds.

\end{itemize}
\end{proposition}

It should be mentioned that (b) can also be derived
from Naftalevi\v c's result \cite[Theorem 3]{Na56}.

From Corollary \ref{maxfct} we can deduce as well a sufficient condition.

\begin{corollary}\label{CSgeom}
Let $\Lambda\subset\D$ be Blaschke. 
If 
\bea\label{CS}
 \sum_{\lambda\in\Lambda} (1-|\lambda|)\log\frac 1{|B_{\lambda}(\lambda)|}<\infty ,
\eea
then $\Lambda\in \Int N^+$ (and so also $\Lambda\in \Int N$).  
\end{corollary}

\begin{proof} 
Set $a_\lambda=\log(1/|B_{\lambda}(\lambda)|)$ and $u=\sum_\lambda
a_\lambda \chi_{I_{\lambda_\lambda}}$.  By assumption  $u\in L^1(\T)$ and obviously
$M_{\Lambda}\leq u$, hence the result follows from Corollary 
\ref{maxfct}(b). 
\end{proof}

It turns out that  for a certain type of sequences condition \eqref{CS} is
both necessary and sufficient for free interpolation in $N$ and $N^+$  (see
Remark \ref{intns}), so that there is no intrinsic analogue of Corollary
\ref{CSgeom} for $N$.  In this direction we state here the following result.

\begin{proposition}\label{CSgeom-best}
For every sequence of positive numbers $(\eps_n)_n\notin\ell^1$
there exists a Blaschke sequence $\Lambda\notin \Int N$ such that 
$\eps_n\simeq (1-|\lambda_n|)\log  1/|B_{\lambda_n}(\lambda_n)|$
 for all $\lambda_n\in\Lambda$.
\end{proposition}

In the comparison of the different geometric conditions we exploit to some
extent the two extremal cases mentioned previously: $\Lambda$  radial (or in a
finite union of Stolz angles), in Proposition \ref{thmnafta}, and $\Lambda$
``very" tangential  (in Proposition~\ref{CSgeom-best}).

The paper is organized as follows. In Section \ref{condsuff} we prove the
sufficiency of the analytic conditions of Theorems~\ref{thmCNS} and \ref{thmCNSN}.
We essentially use a result by Garnett allowing interpolation by $\Hi$ functions on
sequences which are denser than Carleson sequences, under some decrease assumptions
on the interpolated values. 

In Section~\ref{NECESSITY} we study the necessity part of Theorems \ref{thmCNS}
and \ref{thmCNSN}.  We first observe that in the product  $B_{\lambda}$ appearing in
(\ref{CNS}), only the factors $b_{\lambda}(\lambda')$ with $\lambda'$ close to
$\lambda$ are relevant. Then we split the sequence into four pieces,  thereby
reducing the interpolation problem, in a way, to that on separated sequences.

The trace space characterisation will be discussed in  Section \ref{trace}.

Section \ref{GeomCond} is devoted to the proofs of
Propositions \ref{thmnafta} and \ref{CSgeom-best}.

In the final section, we exploit the reasoning of Section \ref{condsuff} to construct
non-Carleson interpolating sequences for ``big'' Hardy-Orlicz classes.

{\bf Acknowledgements}.
We would like to thank Pascal Thomas for some stimulating discussions.

\section{Proof of the sufficient conditions}\label{condsuff}

For a given Blaschke sequence $\Lambda\subset \D$ set
$\delta_\lambda=|B_{\lambda}(\lambda)|$. The key result in the proof of the
sufficient condition is the following theorem by Garnett \cite{Gar77}, that we
cite for our purpose in a slightly weaker form (see also \cite{Nik02} as a
general source, in particular C.3.3.3(g) (Volume 2) for more results of this
kind). 

\begin{theorem*}
Let $\vp:[0,\infty)\longrightarrow [0,\infty)$ be a decreasing function such that
$\int_0^{\infty}\vp(t)\,dt<\infty$. If a sequence $(a_\lambda)_\lambda$ satisfies
\beqa
 |a_\lambda|\le \delta_\lambda\vp (\log\frac{e}{\delta_\lambda}),\quad \lambda
 \in \Lambda,
\eeqa
then there exists a function $f\in \Hi$ such that
$f(\lambda)=a_\lambda$ for all $\lambda\in\Lambda$.
\end{theorem*}

As we have already noted in Remark~\ref{l-infinit}, in order to have free
interpolation in the Nevanlinna and Smirnov classes, it is sufficient that
$\loi\subset N|\Lambda$ and $\loi\subset \NL$ respectively.  Our aim will be
to accommodate the decrease given in Garnett's result by an appropriate
function in $N$ or $N^+$. This is the crucial step in the proof of  the
sufficient conditions of Theorems~\ref{thmCNS} and \ref{thmCNSN}, and it  
occupies its main part.

\begin{proof}[Proof of sufficiency of \eqref{CNS} and \eqref{CNSN}]  
The proof will be presented for the more difficult case of the
Nevanlinna class. So, assume that $g$ and $S$ are as in \eqref{CNSN}.

As in the comment after Theorem \ref{thmCNSN},
we can assume that the function $w\in L^1(\T)$ associated with $g$ by
\beqa
 g(z)=\exp\left( \int_\T\frac{\zeta+z}{\zeta-z}w(\zeta)\,d\sigma(\zeta)\right),
\eeqa
is negative. Let $\mu_S$ be the singular measure associated with $S$. 
The measure $d\mu=-w\, dm+d\mu_S$ is finite and postive, hence
the function
\beqa
 h(z)= \int_\T\frac{\zeta+z}{\zeta-z}d\mu(\zeta)
\eeqa
is holomorphic  with positive real part
in $\D$ (in fact $h=-\log (gS)$). By Smirnov's theorem, $h$ is an outer
function in some $H^p$, $p<1$, and therefore in $N^+$ (see  \cite{Nik02}, in
particular A.4.2.3 (Volume 1)).  
By assumption we have $\log (1/\delta_{\lambda})\le \log (1/|g(\lambda)S(\lambda)|)
=\textrm{Re}\; h(\lambda)$, $\lambda\in\D$.

Take now $\vp(t)=(1+t)^{-2}$, which obviously satisfies the hypothesis of
Garnett's theorem, and set $H=(2+h)^2$, which is still outer in $N^+$. We have the
estimate
\beqa
 |H(\lambda)| = |2+h(\lambda)|^2
  \ge (2+\textrm{Re}\; h(\lambda))^2
  \ge (1+\log \frac e{\delta_\lambda})^2
 =\frac{1}{\vp(\log(e/\delta_\lambda))},
\eeqa
hence the sequence $(\gamma_\lambda)_\lambda$ defined by
\beqa
 \gamma_\lambda=\frac{1}{H(\lambda)\vp(\log(e/\delta_\lambda))}, 
 \quad \lambda\in\Lambda,
\eeqa
is bounded by $1$. 

In order to interpolate $\omega=(\omega_\lambda)_\lambda\in\loi$ by a
function in $N$, split
\beqa
 \omega_\lambda= 
 \big(\omega_\lambda \gamma_\lambda\frac{g(\lambda)S(\lambda)}{\delta_\lambda} 
 \delta_\lambda\vp(\log\frac  e{\delta_\lambda})\big) \cdot
  \frac{H(\lambda)}{g(\lambda)S(\lambda)}.
\eeqa

Since by hypothesis $(\omega_\lambda\gamma_\lambda g(\lambda) 
S(\lambda)/\delta_\lambda)_\lambda$ is bounded, we
can apply Garnett's result to interpolate the sequence 
\beqa
 a_\lambda=\omega_\lambda\gamma_\lambda 
 \frac{g(\lambda)S(\lambda)}{\delta_\lambda} \delta_\lambda 
 \vp(\log\frac e{\delta_\lambda}),\quad \lambda\in \Lambda,
\eeqa
by a function $f\in\Hi$. Now $F=fH/gS$ is a function in $N$ 
with $F|\Lambda=\omega$.

The proof for the Smirnov case is obtained by deleting all the appearences of
$S(\lambda)$ and the singular measure $\mu_S$.
\end{proof}

\section{Proof of the necessary conditions}\label{NECESSITY}

We first show that in order to construct the appropriate function estimating
$|B_{\lambda}(\lambda)|$ from below we only need to consider the factors given by
points $\mu\in\Lambda$ which are close to $\lambda$. This is in accordance with the
results for some related spaces of functions \cite[Theorem 1]{HM2}.

\begin{proposition}\label{propsep}
Let $\Lambda$ be a Blaschke sequence. There exists 
an outer function $g\in N^+$ such that
\beqa
 \prod_{{\mu:|b_{\lambda}(\mu)|\ge 1/2}}|b_{\lambda}(\mu)|
 \ge |g(\lambda)|,\quad
 \lambda\in\Lambda.
\eeqa
\end{proposition}

It is clear from the proof that the constant $1/2$ can be replaced by
any  $\delta\in  (0,1)$. Of course this implies
Corollary~\ref{separated}.

\begin{proof}
Consider the intervals $I_{\lambda}$ defined in \eqref{shadow}.
By the Blaschke condition, the function 
\beqa
 w(\zeta)=-\sum_{\lambda\in\Lambda}\chi_{I_{\lambda}}(\zeta)
\eeqa
belongs to $ L^1(\T)$. Thus the function
\beqa
 g(z)=\exp\left( c \int_\T \frac{\zeta+z}{\zeta-z} w(\zeta)
 \,d\sigma(\zeta)\right),\qquad c>0,
\eeqa
is outer in $N^+$, and it is sufficient to prove that
\beqa
 cP[w](\lambda)=-c\sum_{\mu\in\Lambda}
 \int_{I_\mu}\frac{1-|\lambda|^2}{|\zeta-\lambda|^2}d\sigma(\zeta)
  \le \log \prod_{|b_{\lambda}(\mu)|\ge 1/2}|b_{\lambda}(\mu)|
   =  - \sum_{|b_{\lambda}(\mu)|\ge 1/2} \log \frac{1}{|b_{\lambda}(\mu)|}
\eeqa   
for some constant $c>0$.
   
Using that $\log \frac{1}{t}\simeq 1-t^2$ for $t\in [1/2,1]$, and $
1-|b_{\lambda}(\mu)|^2= \frac{(1-|\lambda|^2)(1-|\mu|^2)}
{|1-\overline{\lambda}\mu|^2}$ we see that
it suffices to prove that
\beqa
 \frac{(1-|\lambda|^2)(1-|\mu|^2)}
   {|1-\overline{\lambda}\mu|^2} 
  \leq c \int_{I_{\mu}} \frac{1-|\lambda|^2}{|\zeta-\lambda|^2}\,
  d\sigma(\zeta) .
\eeqa
We consider two situations. Let $\mu^*=\mu/|\mu|$ and define
\[
D(\mu^*)=\{z\in\D : |z-\mu^*|\leq 2(1-|\mu|)\}.
\]
If $\lambda\in D(\mu^*)$ and $\zeta\in I_{\mu}$  
then $|\zeta-\lambda|\leq 3(1-|\mu|)$, and so 
\[
\int_{I_{\mu}} \frac{d\sigma(\zeta)}{|\zeta-\lambda|^2}\geq \frac 13
\int_{I_{\mu}} \frac{d\sigma(\zeta)}{(1-|\mu|)^2}\simeq\frac 1{1-|\mu|}\gtrsim 
\frac{1-|\mu|}{|1-\overline{\lambda}\mu|^2}.
\]
If $\lambda\notin D(\mu^*)$ and $\zeta\in I_{\mu}$  
then $|\zeta-\lambda|\simeq |\mu-\lambda| $, and so 
\[
\int_{I_{\mu}} \frac{d\sigma(\zeta)}{|\zeta-\lambda|^2}\simeq 
\int_{I_{\mu}} \frac{d\sigma(\zeta)}{|\lambda-\mu|^2}
\simeq\frac {1-|\mu|}{|\lambda-\mu|^2}\gtrsim  
\frac{1-|\mu|}{|1-\overline{\lambda}\mu|^2}.
\]
\end{proof}

\begin{proof}[Proof of the necessity of \eqref{CNS} and \eqref{CNSN}]
We split the sequence into four pieces: $\Lambda=
\bigcup_{i=1}^4 \Lambda_i$ such that each piece $\Lambda_i$ lies in a union
of dyadic ``squares'' which are uniformly separated from each other.
More precisely, consider, for $n\in\N$ and $ 0\le k \le 2^n-1$
the dyadic ``squares":
\beqa
 Q_{n,k}=\bigl\{r\zeta\in \D:1-\frac{1}{2^{n}}\le r < 1-\frac{1}{2^{n+1}},
   \zeta\in \frac{2\pi}{2^n}[k,k+1)\bigr\}.
\eeqa
Set
\beqa
 \Lambda_1=\bigcup_j (\Lambda\cap Q^{(j)}),
\eeqa
where the family $\{Q^{(j)}\}_j$ is given by $\{Q_{2n,2k}\}_{n,k}$
(for the remaining three sequences we respectively choose
$\{Q_{2n ,2k+1}\}_{n,k}$, 
$\{Q_{2n+1,2k}\}_{n,k}$ and
$\{Q_{2n+1,2k+1}\}_{n,k}$).
In order to avoid technical difficulties we count only those
$Q^{(j)}$ containing points of $\Lambda$ (in case $\Lambda_j$
is empty there is nothing to prove).
In what follows we will argue on one sequence, say $\Lambda_1$. 
The arguments are the same for the other sequences.

Our first observation is that, by construction,
\beqa
 \rho(Q^{(j)},Q^{(l)}):=\inf_{z\in  Q^{(j)}, w\in Q^{(l)}}|b_z(w)|
 \ge\delta>0,\quad j\neq l,
\eeqa
for some fixed $\delta$. Also, the closed rectangles 
$\overline{Q^{(j)}}$ are compact in
$\D$ so that $\Lambda_1\cap Q^{(j)}\subset \Lambda\cap Q^{(j)}$ can only
contain a finite number of points (they contain at least one point,
by assumption). 
Therefore
\beqa
 0<m_j :=\min_{\lambda\in \Lambda_1\cap Q^{(j)}} |B_{\lambda}(\lambda)|
\eeqa
(note that we consider the entire Blaschke product $B_\lambda$
associated with $\Lambda\setminus\{\lambda\}$). Take
$\lambda_j^{1}\in Q^{(j)}$ such that
$m_j=|B_{\lambda_j^{1}}(\lambda_j^{1})|$.

Assume now that $\Lambda\in \Int N$.
Since $\loi\subset N\vert\Lambda$, there exists a function $f_1\in N$
such that 
\beqa
 f_1(\lambda)=
 \begin{cases}
1\quad &\text{if}\quad \lambda\in \{\lambda_j^{1}\}_j\\
0\quad &\text{if}\quad \lambda\notin \{\lambda_j^{1}\}_j.
\end{cases}
\eeqa
By the Riesz-Smirnov factorisation we have
\beqa
 f_1=B_{\Lambda\setminus \{\lambda_j^{1}\}_j}\frac{h_1}{h_2 T_2},
\eeqa
where $T_2$ is singular inner, $h_1$ is some function in $\Hi$ and $h_2$ is
outer in $\Hi$. 
Again, we can assume
$\|h_i\|_{\infty}\le 1$, $i=1,2$. Hence
\beqa
 1 = f_1(\lambda_k^{1})
  \le |B_{\Lambda\setminus \{\lambda_j^{1}\}_j}(\lambda_k^{1})|\cdot
   \frac{1}{|h_2(\lambda_k^{1}) T_2(\lambda_k^{1})|},
\eeqa
and
\beqa
 |B_{\Lambda\setminus \{\lambda_j^{1}\}_j}(\lambda_k^{1})|
  \ge {|h_2(\lambda_k^{1}) T_2(\lambda_k^{1})|},\quad  k\in\N.
\eeqa

Since $h_2 T_2$ does not vanish and is bounded above by 1, the function
$\log |h_2 T_2|$ is a negative harmonic function. By Harnack's
inequality, there exists an absolute constant $c\ge 1$ such that
\beqa
 \frac{1}{c} |\log |h_2(\lambda_k^{1}) T_2(\lambda_k^{1})||
 \le |\log |h_2(z) T_2(z)|| \le c |\log|h_2(\lambda_k^{1}) T_2(\lambda_k^{1})||,
 \quad z\in Q^{(k)}, 
\eeqa
hence
\beqa
|h_2(\lambda_k^{1})T_2(\lambda_k^{1})|^{c}\le |h_2(z) T_2(z)|
 \le |h_2(\lambda_k^{1}) T_2(\lambda_k^{1})|^{1/c},\quad z\in Q^{(k)}.
\eeqa 
This yields
\bea\label{estimate4}
 |(h_2 T_2)^c(\mu)|\le |(h_2 T_2)(\lambda_k^{1})|\le 
 |B_{\Lambda\setminus \{\lambda_j^{1}\}_j}(\lambda_k^{1})|
\eea
for every $\mu\in \Lambda_1\cap Q^{(k)}$.

Let us now exploit Proposition \ref{propsep}. By construction,
the sequence $\{\lambda_j^{1}\}_j\subset\Lambda_1$ is 
separated. Therefore, there
exists an outer function $G_1$ in the Smirnov class such that
\beqa
 |B_{\{\lambda_j^{1}\}_j\setminus \{\lambda_k^{1}\}}(\lambda_k^{1})|
 \ge |G_1(\lambda_k^{1})|,\quad k\in\N.
\eeqa
Again, $G_1$ is a quotient of two bounded outer functions and 
we can suppose that $G_1$ is outer in $\Hi$ with
$\|G_1\|_{\infty}\le 1$. Also, we can use Harnack's inequality as above
to get
\beqa
 |G_1(\lambda_k^{1})|\ge |G_1^c(\mu)|
\eeqa
for every $\mu\in \Lambda_1\cup Q^{(k)}$.
This together with (\ref{estimate4}) and our definition
of $\lambda_k^{1}$ give
\beqa
 |B_{\Lambda\setminus\{\mu\}}(\mu)|
 &\ge& |B_{\Lambda\setminus\{\lambda_k^{1}\}}(\lambda_k^{1})|
 =|B_{\Lambda\setminus \{\lambda_j^{1}\}_j}(\lambda_k^{1})|
 \cdot |B_{\{\lambda_j^{1}\}_j\setminus \{\lambda_k^{1}\}}
   (\lambda_k^{1})|\\
 &\ge& |(h_2 T_2)^c(\mu)|\cdot |G_1^c(\mu)|
\eeqa
for every $\mu\in Q^{(k)}$ and $k\in\N$.
Setting now $g_1=(h_2G_1)^c$ and $S_1=T_2^c$ we get the claim for all points
in $\Lambda_1$. Note also that, by construction, $g_1$ is outer with
$\|g_1\|_{\infty}\le 1$ and $S_1$ is singular inner.

To finish the proof, construct in a similar way functions $g_i$, $S_i$
for the sequences $\Lambda_i$, $i=2,3,4$ and define the products
\beqa
 g=\prod_{i=1}^4 g_i\quad\textrm{and}\quad S=\prod_{i=1}^4 S_i.
\eeqa
Of course $g$ is outer in $\Hi$, and $S$ is singular
inner. So, whenever $\mu\in\Lambda$, there exists
$k\in \{1,2,3,4\}$ such that $\mu\in\Lambda_k$, and hence
\beqa
 |B_{\lambda}(\lambda)|\ge |g_k(\lambda)S_k(\lambda)|\ge |g(\lambda)S(\lambda)|.
\eeqa
The proof for $N^+$ follows the same lines, just disregarding the singular inner
factors.
\end{proof}

\section{The trace space}\label{trace}

In this short section we  prove the trace space characterisation of free
interpolation given in Theorems \ref{thmCNS} and \ref{thmCNSN}. 

In order to see that (d) of both theorems implies free interpolation it suffices to
observe that  $\ell^{\infty}\subset l_{N^+}\subset  l_N$ and to use Remark
\ref{l-infinit}.

For the proof of the converse, we will only consider the situation in the Nevanlinna
class, since the case of the Smirnov class is again obtained by removing the
singular part of the measure and the singular inner factors. 

Assume that $(a_\lambda)_\lambda\in N|\Lambda$ and that $f\in N$ is such that
$f(\lambda)=a_\lambda$, $\lambda\in\Lambda$. Since $f$ can be written as
$f=f_1/S_2f_2$, where $f_1\in \Hi$, $\|f_1\|_{\infty}\le 1$, $S_2$ is singular inner
with associated singular measure $\mu_S$, and $f_2\in \Hi$ is an outer function with
$\|f_2\|_{\infty}\le 1$, we can define the positive finite measure $\mu=\log
(1/|f_2|)\, dm+d\mu_S$, that obviously satisfies
$P[\mu](\lambda)\geq\log^+|a_\lambda|$, $\lambda\in\Lambda$.

Conversely, suppose that $(a_{\lambda})_{\lambda}$ is
such that there is a positive finite measure $\mu$ with $P[\mu](\lambda)\ge
\log^+|a_{\lambda}|$. The Radon-Nikodym decomposition of $\mu$ is given by
$d\mu=w\,dm+d\mu_S$, where $w\in L^1(\T)$ is positive and $\mu_S$ is a positive
finite singular measure. Let $S$ be the singular inner function associated with
$\mu_S$, and let $f$ be the function defined by
\beqa
 f(z)=\exp\left(\int_{\T}\frac{\zeta+z}{\zeta-z} w(\zeta)\,d\sigma(\zeta)\right),
\quad z\in\D.
\eeqa
By definition, $f$ is outer in $N^+$ and $F=f/S\in N$. Clearly,
$\log^+|a_{\lambda}|\le \log |F(\lambda)|$, thus $|a_{\lambda}|\le |F(\lambda)|$. 
Since $N|\Lambda$ is ideal by assumption, 
there exists  $f_0\in N$ interpolating $(a_{\lambda})_{\lambda}$. 
\hfill \qedsymbol

\section{Proofs of Propositions \ref{thmnafta} and \ref{CSgeom-best}}\label{GeomCond}

\begin{proof}[Proof of Proposition \ref{thmnafta}]

It is enough to consider the case where $\Lambda$ is contained in only one Stolz
angle. Indeed, if $\Lambda=\bigcup_{i=1}^n\Lambda_i$ with $\Lambda_l\subset  
\Gamma_{\zeta_l}$, $l=1,\ldots,n$, and $\zeta_i\neq\zeta_j$, then $\lim\limits_{z\to
\zeta_i, z\in \Gamma_{\zeta_i}} |B_{\Lambda_j}(z)|=1$, $j\neq i$, so that $\log
1/|B_{\lambda}(\lambda)|$ behaves asymptotically  like $\log
1/|B_{\Lambda_i\setminus\lambda}(\lambda)|$
in $\Gamma_{\zeta_i}$ (here $\lambda\in \Lambda_i$).

Also, we can assume that the sequence is radial
(this means that we replace the initial sequence by 
one which is in a uniform pseudo-hyperbolic neighbourhood of the 
initial one; by Harnack's inequality such a perturbation does not 
change substantially the behaviour of positive harmonic functions).

By Corollary \ref{CNgeom} it remains to prove the sufficiency of the conditions.
Let us first show that
condition \eqref{CN} implies interpolation in
$N^+$ . In order to construct a function $w\in L^1(\T)$ meeting the requirement of Theorem
\ref{thmCNS}(c) assume that $\Lambda=\{\lambda_n\}_n\subset [0,1)$ is arranged in
increasing order and set $\tilde{\eps}_n=(1-|\lambda_n|)\log (1/|B_{\lambda_n}
(\lambda_n)|)$. Clearly there exists a decreasing sequence $(\eps_n)_n$
with $\tilde{\eps}_n\le\eps_n$, $n\in \N$, and $\lim_n \eps_n=0$. Now, if
$I_n=I_{\lambda_n}$ are the arcs defined in \eqref{shadow}, $J_n=I_n\setminus
I_{n+1}$ and $\beta_n= \eps_n-\eps_{n+1}$, set  
\beqa
 w(\zeta)=\sum_n \frac{\beta_n}{|J_n|} \chi_{J_n}(\zeta),\quad
 \zeta\in\T,
\eeqa
where $|J_n|$ denotes the Lebesgue measure of the set $J_n$.
Then $w\in L^1(\T)$, and 
\begin{eqnarray*}
 P[w](\lambda_n) &\ge&  \int_{I_n} P(\lambda_n,\zeta)
   \sum_k \frac{\beta_k}{|J_k|} 
  \chi_{J_k}(\zeta) \,d\sigma(\zeta)
 \gtrsim \sum_{k\ge n}\frac{\beta_k}{|J_k|}
   \frac{1}{(1-|\lambda_n|)}\int_{J_k}\,d\sigma(\zeta)\\
 &=&\frac{\sum_{k\ge n} \beta_k}{(1-|\lambda_n|)}
 =\frac{\eps_n}{1-|\lambda_n|}
  \ge \frac{\tilde{\eps}_n}{1-|\lambda_n|}
=\log \frac 1 {|B_{\lambda_n} (\lambda_n)|}.
\end{eqnarray*}
This and Theorem \ref{thmCNS} prove the assertion.

The proof for the Nevanlinna class is even simpler. 
Set $d\mu_s=  \delta_{1}$, the Dirac mass on $1\in\T$. From  \eqref{CNN} we get 
\[
 \log \frac 1{|B_{\lambda_n}(\lambda_n)|}\lesssim \frac{1}{1-|\lambda_n|}
 \lesssim  P[\mu_s](\lambda),
\]
and we finish by applying Theorem \ref{thmCNSN}.
\end{proof}

Similar ideas are used in the proof of Proposition \ref{CSgeom-best},
which will be a consequence of the following auxiliary result. For
$\lambda\in \D$, denote
\beqa
 K_{\lambda}=\{\zeta\in\T:|\arg\zeta-\arg\lambda|\le \pi\sqrt{1-|\lambda|}\}.
\eeqa 

\begin{proposition}\label{noouter}
Let $\Lambda\subset \D$ be such that 
\beqa 
 K_\lambda\cap K_{\lambda'} =\emptyset,\quad \lambda, \lambda'\in\Lambda,
 \quad \lambda\neq \lambda'.
\eeqa
Then,  
for no sequence of positive numbers $(\eps_\lambda)_\lambda$ with
$\sum_\lambda \eps_\lambda=\infty$ there exists a finite positive measure
$\mu$ such that
$\eps_\lambda/(1-|\lambda|)\le P[\mu](\lambda)$ for all $\lambda\in\Lambda$
\end{proposition}

Notice that the hypothesis implies automatically
that $\Lambda$ is a Carleson sequence.

\begin{proof}
Suppose on the contrary that there exists a finite positive measure $\mu$
satisfying the above estimate. Note that
\beqa
 P[\mu](\lambda)
 = \int_{K_\lambda} P(\lambda,\zeta)\,d\mu(\zeta)+
  \int_{\T\setminus K_\lambda} P(\lambda,\zeta)\,d\mu(\zeta)
  \le \int_{K_\lambda} P(\lambda,\zeta)\,d\mu(\zeta)+c_{\mu},
\eeqa
where $c_{\mu}=\mu(\T)$.
So
\beqa
 \eps_\lambda \le c_{\mu} (1-|\lambda|)+
     (1-|\lambda|)\int_{K_\lambda}P(\lambda,\zeta)\,d\mu(\zeta) 
  \le c_{\mu} (1-|\lambda|)+2 \int_{K_\lambda}\,d\mu(\zeta) ,
\eeqa
and consequently
\beqa
 \sum_{\lambda\in\Lambda} \eps_\lambda\le c_{\mu} \sum_{\lambda\in\Lambda} 
 (1-|\lambda|)
 +2 \int_{\bigcup_\lambda K_\lambda}\,d\mu(\zeta) 
 \le c_{\mu} \sum_{\lambda\in\Lambda} (1-|\lambda|)
 +2 c_{\mu}.
\eeqa
Since $\Lambda$ is a Blaschke sequence, this is not possible for
any finite measure $\mu$.
\end{proof}

\begin{proof}[Proof of Proposition \ref{CSgeom-best}]
Let $(\eps_n)_n\not\in \ell^1$.  There is no restriction in assuming   $\lim_n
e^{-\eps_n/(1-|\lambda_n|)}=0$.

Take now $\Lambda_1$ as in the previous proposition. In the neighbourhood of each
$\lambda_n=r_ne^{i\vp_n}\in\Lambda_1$ we fix a point
$\lambda_n^\prime=r_n'e^{i\vp_n}$ such that $r_n'\ge r_n$ and 
$|b_{\lambda_n^\prime}(\lambda_n)|=e^{-\eps_n/(1-|\lambda_n|)}$ (so
$\lambda_n^\prime$ lies on the shortest  segment connecting $\lambda_n$ to $\T$).
Since  $e^{-\eps_n/(1-|\lambda_n|)}$ decreases to 0, by a standard perturbation
argument, the sequence $\Lambda_2=\{\lambda_n^\prime\}_n$ is also Carleson.

Define $\Lambda=\Lambda_1\cup\Lambda_2$. By construction
\beqa
 |B_{\lambda}(\lambda)|\simeq e^{-\eps_{\lambda}/(1-|\lambda|)},
 \quad \lambda\in\Lambda,
\eeqa
where $\eps_{\lambda}=\eps_n$ both for $\lambda=\lambda_n$ and 
$\lambda=\lambda^\prime_n$.

The previous proposition and Theorem \ref{thmCNSN}(c) show that
$\Lambda\not\in\Int N$.
\end{proof}

\begin{remark}\label{intns}
Let $\Lambda=\{\lambda_{n,k}\}_{n\in\N, k\le N}$ be a sequence
constructed similarly to the prove above,
i.e. $N\in\N$ is some fixed number,
$\{\lambda_{n,1}\}_n$ satisfies the hypothesis of Proposition~\ref{noouter} and
$|b_{\lambda_{n,k}}(\lambda_{n,l})|\le \delta$, $k,l=1,\ldots,N$, for some fixed 
$\delta\in (0,1)$. Then the previous results show that the
following assertions are equivalent
\begin{itemize}
\item[(i)] $\sum_{\lambda\in\Lambda}(1-|\lambda|)\log 1/|B_{\lambda}(\lambda)|
<\infty$,
\item[(ii)] $\Lambda\in \Int N$,
\item[(iii)] $\Lambda\in \Int N^+$.
\end{itemize}
\end{remark}

\section{Hardy-Orlicz classes}\label{nevanlinna}

Let $\vp:{\R}\lra [0,\infty)$ be a convex, nondecreasing
function satisfying 
\begin{itemize}
\item[(i)] $\lim_{t\to \infty} \vp (t)/t =\infty$
\item[(ii)] $\Delta_2$-condition: $\vp (t+2) \le M \vp (t) +K$, 
$t \ge t_0$ for some constants $M,K \ge 0$ and $t_0\in \R$. 
\end{itemize}
Such a function is called strongly convex (see \cite{RR}),
and one can associate with it the corresponding \emph{Hardy-Orlicz
class}
\beqa
 \hf =\{f\in N^+:\int_\T\vp (\log |f(\zeta)|)\,d\sigma(\zeta) 
 <\infty\},
\eeqa 
where $f(\zeta)$ is the non-tangential boundary value of $f$ at
$\zeta\in {\T}$, which exists almost everywhere. 
In \cite{Har}, the following result was proved.

\begin{theorem*}
Let $\vp$ be a strongly convex function satisfying (i), (ii) and 
the $V_2$-condition: 
\beqa
 2\vp(t)\le \vp(t+\alpha), \quad t\ge t_1
\eeqa 
where $\alpha>0$ is a suitable constant and $t_1\in\R$. Then $\Lambda
\subset {\D}$ is free interpolating for $\hf$ 
if and only if
$\Lambda$ is a Carleson sequence, and in this case 
\beqa
 \hf\vert_{\Lambda}=\{a=(a_\lambda)_\lambda: |a|_{\varphi}=
 \sum_{\lambda} (1-|\lambda|) 
 \varphi(\log |a_\lambda|)<\infty\}.
\eeqa
\end{theorem*}

The conditions on $\vp$ imply that for all
$\hf$ there exist $p,q\in (0,\infty)$ such that $H^p\subset \hf\subset H^q$.
In particular, the $V_2$-condition
implies the inclusion $H^p\subset \hf$ for some $p>0$. This $V_2$-condition has
a strong topological impact on the spaces. In fact, it guarantees that metric
bounded sets are also bounded in the topology of the space (and so the
functional analysis tools still apply in this situation; see
\cite{Har} for more on this and for further references).
It was not clear whether this was only a technical problem or if there
existed a critical growth for $\vp$ (below exponential growth 
$\vp(t)=e^{pt}$ corresponding to $H^p$ spaces) giving a breakpoint in
the behaviour of interpolating sequences for $\hf$. 

We can now affirm that this behaviour in fact changes between exponential
and polynomial growth. Let $\vp$ be a strongly convex function
with associated Hardy-Orlicz space $\hf$. Assume moreover that 
$\vp$ satisfies 
\bea\label{deltaineq}
 \vp(a+b)\le c(\vp(a)+\vp(b)),
\eea
for some fixed constant $c\ge 1$ and for all $a,b\ge t_0$.
The standard example in this setting is $\vp_p(t)=t^p$ for $p>1$.
We have the following result.

\begin{theorem}
Let $\vp:\R\lra [0,\infty)$ be a strongly convex function such that
\eqref{deltaineq} holds. If there exists a positive weight $w\in L^1(\T)$ such that
$\vp \circ w\in L^1(\T)$ and \eqref{CNSmu} holds, then $\Lambda\in \Int \hf$.
\end{theorem}

\begin{proof}
Note first that \eqref{deltaineq} implies that $\hf$
is an algebra contained in $N^+$, hence it is sufficient to interpolate
bounded sequences (see Remark~\ref{l-infinit}). 
By assumption, the outer function defined by
\beqa
 g(z)=\exp \left\{\int_\T \frac{\zeta+z}{\zeta-z}  
  (-w(\zeta))d\sigma(\zeta) \right\}
\eeqa
satisfies $|B_{\lambda}(\lambda)|\ge |g(\lambda)|$, $\liL$.
The reasoning carried out in 
Section \ref{condsuff}
leads to an interpolating function of the form $fH/g$, with 
$f\in \Hi$, and $H=(2+h)^2$ is  
outer in $H^p$ for all $p<1$ (note that the measure $\mu$ defining $h$ here
is absolutely continuous,
in fact $\mu=-w\,dm$). Also, $H^p\subset \hf$ for any $p>0$ by our conditions
on $\vp$. By construction, 
$\int \vp (\log |1/g|)=\int \vp \circ w <\infty$ so that $1/g\in \hf$.
Since $\hf$ is an algebra, we deduce that $fH/g\in\hf$.
\end{proof}

\begin{example}
We give an example of an interpolating sequence for
$\hf$ which is not Carleson,
thus justifying our claim that there
is a breakpoint between Hardy-Orlicz spaces verifying the $V_2$-condition
and those that do not.

Consider the functions $\vp_p$ and let $\Lambda_0=\{\lambda_n\}_n\subset \D$ be
a Carleson sequence  verifying $I_n\cap I_k=\emptyset$, $n\neq k$. Since
$\sum_n (1-|\lambda_n|)<\infty$, there exists a strictly increasing sequence of
postive numbers $(\gamma_n)_n$ such that  $\sum_n
(1-|\lambda_n|)\gamma_n<\infty$ and $\lim_{n\to\infty}\gamma_n=\infty$. 
Setting  $\eps_n=(1-|\lambda_n|)\gamma_n^{1/p}$ and 
\beqa
 u=\sum_n\frac{\eps_n}{1-|\lambda_n|}\chi_{I_n},
\eeqa
we obtain $\int \vp \circ u =\sum_n\frac{\eps_n^p}{(1-|\lambda_n|)^{p-1}}
=\sum_n(1-|\lambda_n|)\gamma_n <\infty$.  As in the proof of Proposition
\ref{CSgeom-best}, we attach a second Carleson sequence $\Lambda_1=\{\lambda^\prime_n\}_n$
such that  the pseudo-hyperbolic distance between corresponding points
satisfies $|b_{\lambda^\prime_n} (\lambda_n)|=e^{-\eps_n/(1-|\lambda_n|)}$.  Since
$\gamma_n\to\infty$ we have $\eps_n/(1-|\lambda_n|)=\gamma_n^{1/p}\to \infty$,
i.e.\ the elements of the sequence $\Lambda=\Lambda_1\cup\Lambda_2$  are
arbitrarily close and $\Lambda$ cannot be a Carleson sequence.  By construction,
condition \eqref{CNSmu} holds (as before, we may possibly have to multiply
$u$ with some constant $c$ to have that condition also in the points
$\lambda^\prime_n$, but this operation conserves the integrability condition), and
therefore $\Lambda\in\Int \hf$. 
\end{example}

\providecommand{\bysame}{\leavevmode\hbox to3em{\hrulefill}\thinspace}
\providecommand{\MR}{\relax\ifhmode\unskip\space\fi MR }


\begin{thebibliography}{BRSHZE}


\bibitem[Gar77]{Gar77} 
{\it J.B. Garnett}, 
{Two remarks on interpolation by bounded analytic functions},
Banach spaces of analytic functions (Proc. Pelczynski Conf., Kent State 
Univ., Kent, Ohio, 1976), pp. 32--40. Lecture Notes in Math., Vol. 604, 
Springer, Berlin, 1977. 

\bibitem[Gar81]{Gar} 
{\it J.B. Garnett}, 
{Bounded analytic functions}, Academic Press, New York, 1981. 

\bibitem[Har99]{Har}
{\it A. Hartmann},  
{Free interpolation in Hardy-Orlicz spaces},
Studia Math. \textbf{135} (1999), no. 2, 179--190. 

\bibitem[HaMa01]{HM2}
{\it A. Hartmann \& X. Massaneda},
{Interpolating sequences for holomorphic 
functions of restricted growth}, to appear in Ill. J. Math.

\bibitem[McC92]{McC92}
{\it J. McCarthy},
{Topologies on the Smirnov class},
J. Funct. Anal. \textbf{104} (1992), no. 1, 229--241. 

\bibitem[Na56]{Na56} 
{\it A.G. Naftalevi\v c},
{On interpolation by functions of
bounded characteristic (Russian)},  Vilniaus Valst. Univ. Moksl\c u Darbai.
Mat. Fiz. Chem. Moksl\c u \textbf{Ser. 5} (1956), 5--27. 

\bibitem[Nik86]{niktr} 
{\it N.K. Nikolski [Nikol'ski\u{\i}]},
{Treatise on the shift operator}, Springer-Verlag, Berlin etc., 1986. 

\bibitem[Nik02]{Nik02} 
\bysame
\emph{Operators, functions, and systems: an easy reading. Vol. 1,
Hardy, Hankel, and Toeplitz; Vol.2, Model Operators and Systems}, 
Mathematical Surveys and Monographs, 92 and 93. 
American Mathematical Society, Providence, RI, 2002. 

\bibitem[RosRov]{RR}
{\it M. Rosenblum \& J. Rovnyak},
\emph{Hardy classes and operator theory}, 
Oxford Mathematical Monographs. Oxford Science Publications. 
The Clarendon Press, Oxford University
Press, New York, 1985.

\bibitem[ShHSh]{ShHSh}
{\it H.S. Shapiro \& A.L. Shields},
{On some interpolation problems for 
analytic functions}, Amer. J. Math., \textbf{83} (1961), 513--532.

\bibitem[ShSh]{ShSh} 
{\it J. Shapiro \& A. Shields}, 
{Unusual topological
properties of the Nevanlinna class}, Amer. J. Math. \textbf{97} (1975),
915--936. 

\bibitem[Ya74]{yana2} 
{\it N. Yanagihara},
{Interpolation theorems for the
class $N^+$}, Illinois J. Math., \textbf{18} (1974), 427--435.
\end{thebibliography}
\end{document}